\documentclass[11pt]{article}

\title{\bf The Cauchy problem and integrability of a modified Euler-Poisson equation\footnote{to appear in {\em Transactions of the American Mathematical Society }(in press)}}

\author{Feride T\i\u{g}lay}

\usepackage{amsmath, amssymb, amsthm}

\def\eps{\varepsilon}
\def\ga{\gamma}
\def\gi{\gamma^{-1}}

\def\z{\zeta}

\def\diff{\mathcal{D}^{s}}
\def\z{\zeta}

\def\dt{\partial_{t}}
\def\dx{\partial_{x}}

\def\sob{H^{s}}
\def\sobb{H^{s-1}}

\def\tor{\mathbb{T}}

\def\del{\Delta}

\def\gr{\nabla}
\def\reel{\mathbb{R}}
\def\torm{\mathbb{T}^{m}}
\def\lop{\Lambda^{-2}}

\def\di{\partial_{i}}
\def\dj{\partial_{j}}

\def\dm{\partial_{m}}

\def\D{\mathcal{D}}
\def\P{\mathcal{P}}

\def\po{\phi_{1}}
\def\pt{\phi_{2}}
\def\to{\theta_{1}}
\def\tt{\theta_{2}}
\def\fn{\frac{\delta F}{\delta n}}
\def\fv{\frac{\delta F}{\delta v}}
\def\hn{\frac{\delta H}{\delta n}}
\def\hv{\frac{\delta H}{\delta v}}
\def\grad{\mathrm{grad}}
\def\div{\mathrm{div}}

\newtheorem {lema} {Lemma}
\newtheorem {prop}{Proposition}
\def \pf {\textit{Proof. }}

\newtheorem {rema}{Remark}
\newtheorem {teor}{Theorem}

\begin{document}
\maketitle
\begin{abstract}
We prove that the periodic initial value problem for a modified Euler-Poisson equation is well-posed for initial data in $\sob (\torm)$ when $s>m/2+1$. We also study the analytic regularity of this problem and prove a Cauchy-Kowalevski type theorem. After presenting a formal derivation of the equation on the semidirect product space $ \mathrm{Diff} \ltimes C^{\infty}(\tor)$ as a Hamiltonian equation, we concentrate to one space dimension ($m=1$) and show that the equation is bihamiltonian.
\end{abstract}

In this paper we study the periodic Cauchy problem for the modified Euler-Poisson equation (mEP)
\begin{equation*}
\begin{array}{l}
\dt n + \mathrm{div} (nv)=0 \\ 
\dt v + (v.\gr) v +  \mathrm{grad}\lop  n =0, \ x\in\torm, \ t\in\mathbb{R}.
\end{array} 
\tag{mEP}
\end{equation*}
as well as its hamiltonian structure and integrability.

The equation (mEP) is related to the Euler-Poisson equation
\begin{equation}  \begin{array}{l}
\dt n +  \mathrm{div}(nv)=0 \\ 
\dt v + (v.\gr) v + \mathrm{grad}  \phi =0 \\ 
\del \phi - e^{\phi} + n=0, x\in\torm, t\in\mathbb{R}
\end{array} 
\label{eq:iapintro}
\end{equation}
which describes the fluctuations in the ion density of a two-component plasma of positively charged ions and
negatively charged electrons (therefore it is also called ion acoustic plasma equation \cite{LiSat}).

Linearizing the operator $N(\phi)=e^{\phi}-\del\phi$ at $\phi=0$ in the
Euler-Poisson equation (\ref{eq:iapintro}), we obtain the local form of the modified Euler-Poisson equation (mEP):
\begin{equation}  \begin{array}{l}
\dt n +  \mathrm{div}(nv)=0, \\ 
\dt v + (v.\gr) v + \mathrm{grad}  \phi =0, \\ 
\del \phi - \phi + n=0,\  x\in\torm, \ t\in\mathbb{R}.
\end{array}  \label{eq:mlineer} 
\end{equation}
The equation (\ref{eq:mlineer}), like the Euler-Poisson equation (\ref{eq:iapintro}), admits an approximation which preserves dispersion and leads to KdV (see remark \ref{rem:1}). 
Inverting the operator $\Lambda :=I-\Delta$, we write the system in (\ref{eq:mlineer}) in the nonlocal form (mEP). 

Besides its relation to the Euler-Poisson equation (\ref{eq:iapintro}), the modified Euler-Poisson equation is also remarkable for its bihamiltonian structure in one space dimension that we describe here.

Many different techniques have been developed based on Picard's contraction theorem on Banach spaces in the study of nonlinear partial differential equations. One approach originated  in an
observation of V. Arnold \cite{Arn} that the initial value problem for the classical Euler equations of a perfect fluid
can be stated as a problem of finding geodesics on the group of volume preserving
diffeomorphisms. Subsequently, this observation was used by D.G. Ebin and J. Marsden in \cite{EMa} who
developed the necessary functional analytic tools and established sharp local well-posedness results for the Euler
equations in a class of Sobolev spaces. 

The first section of this work is devoted to develop an appropriate analytic framework for the modified Euler-Poisson equation (mEP) using a similar approach and prove the following theorem.

\begin{teor}
For $s>m/2+1$, given any initial data $(n_{0}, v_{0}) \in \sobb(\torm,\mathbb{R}) \times
\sob(\torm,\mathbb{R}^{m})
$, there exists a
$T>0$ and a unique solution
$(n,v)$ to the Cauchy problem for the modified Euler-Poisson equation (mEP) such that 
\[ v \in C([0,T), \sob(\torm, \mathbb{R}^{m})) \cap C^{1}([0,T), \sobb(\torm,\mathbb{R}^{m})) \]
and 
\[ n \in C([0,T), \sobb(\torm,\mathbb{R})) \cap C^{1}([0,T), H^{s-2}(\torm,\mathbb{R})) \] 
and the solution
$(n,v)$ depends continuously on the initial data
$(n_{0},v_{0})$.
\label{th1}
\end{teor}

Another powerful tool in the study of partial differential equations is the Cauchy-Kowalevski theorem. An abstract
version of this theorem was developed by L.V. Ovsjannikov \cite{Ovs1,Ovs2}, F. Treves
\cite{Tre}, L. Nirenberg \cite{Nir}, T. Nishida \cite{Nis} and M.S. Baouendi and C.
Goulaouic
\cite{BG} among others and subsequently applied to the Euler and Navier-Stokes equations. 

The study of analytic
regularity of solutions to the Camassa-Holm equation by A. Himonas and G. Misio\l ek \cite{HM1}, \cite{HM3} using this abstract
theorem led us to investigate the analytic regularity for (mEP). In Section 2, 
we prove the existence and uniqueness of local analytic solutions to the Cauchy problem for the
equation (mEP).

\begin{teor}
If the initial data $(n_{0},v_{0})$ is analytic on $\torm\times\torm$ then there exists
an $\eps>0$ and a unique solution $(n,v)$ of the Cauchy problem for the equation (mEP) that is analytic both in $x$ and $t$
on $\torm\times\torm$ for all $t$ in $(-\eps, \eps)$.
\label{th:iap_analy_reg}
\end{teor}

This result can be viewed as a Cauchy-Kowalevski type result for the equation (mEP). Even though the equation (mEP) admits an approximation by the Korteweg-De Vries equation, the analytic regularity results for the two equations are quite different. In contrast with the Korteweg-De Vries equation whose solutions are analytic in the space variable for all time but not analytic in the time variable (see \cite{Tru}, \cite{KaM}), the solutions to the modified Euler-Poisson equation are analytic in both space and time variables.

In the third section we derive the equation (mEP) as a Hamiltonian equation on the semidirect product space
$\mathrm{Diff}(\torm) \ltimes C^{\infty}(\torm)$ following the treatment of V. Arnold and B. Khesin in \cite{ArKh} and J.
Marsden, T. Ratiu and A. Weinstein in \cite{MRW} of the Hamiltonian formalism related to fluid and gas dynamics.
Then we concentrate to the one space dimension $m=1$ and prove the following theorem.
\begin{teor}
For $m=1$ the modified Euler-Poisson equation (mEP) is bihamiltonian with the pair of hamiltonian functionals
\[ H_{1}=\int \frac{1}{2}(v^{2}n+(\lop\dx n)^{2}+(\lop n)^{2})dx
\]
and
\[ H_{2}=\int nv \ dx.
\]
\label{th:integrable}
\end{teor}
In the proof we use prolongations to check the compatibility of the induced Poisson brackets by these hamiltonian structures. In particular the modified Euler-Poisson equation (mEP) can be derived as a Hamiltonian equation on the semidirect product space of the Virasoro algebra with the smooth functions on the torus $vir\ltimes C^{\infty}(\tor)$ along with a nonlocal hierarchy of equations called Hunter-Zheng equations (see \cite{BDP} for the bihamiltonian structure of the Hunter-Zheng equations).
\begin{rema}
The Korteweg-de Vries equation (KdV) can be derived as an approximation to the Euler-Poisson equation by a perturbation analysis (see \cite{Sattinger}). Using this approach it is straightforward to obtain  an approximation to the system of equations in (\ref{eq:mlineer}) which preserves the dispersion and leads to KdV. 
\label{rem:1}
\end{rema}

\section{Local well-posedness in Sobolev spaces}

In this section we study the Cauchy problem for the modified Euler-Poisson equation (mEP)
where $n=n(t,x):\reel \times \tor^{m}\rightarrow\reel$ and $v=v(t,x):\reel \times
\tor^{m}\rightarrow\mathbb{R}^{m}$ and $\lop=(I-\del)^{-1}$ is the Bessel potential. 

In order to prove Theorem \ref{th1} we use the method 
of first restating the problem as an initial value problem for an ordinary differential equation on
the group of diffeomorphisms of Sobolev class $\sob$ and then applying the existence theorem for vector fields
on Banach manifolds. 
\begin{prop}
For $s>m/2+1$, a pair $(n,v) \in \sobb(\torm,\mathbb{R})\times\sob(\torm,\mathbb{R}^{m})$ is a solution to
the Cauchy problem for (mEP) with initial data $(n_{0},v_{0})$ if and only if  $v=\eta \circ \gi$
and
$n=\z\circ\gi$ where $(\ga,\z,\eta)$ is a solution to
\begin{align}
& \dt \z = -((\z \circ \gi)  \mathrm{div} (\eta \circ \gi)) \circ \ga, \nonumber \\ 
& \dt \eta =-(\mathrm{grad} (\lop (\z \circ \gi))) \circ \ga, \label{eq:yeni2} \\ 
& \dt\ga =\eta, \nonumber  
\end{align}
with initial conditions $\z (0,x)=n_{0}(x) , \eta(0, x)=v_{0}(x) , \ga(0, x)=id_{x}$.  
\label{pro:ode_iap}
\end{prop}

Therefore the Cauchy problem for (mEP) can be reformulated as an initial value problem for the ordinary
differential equation 
\begin{equation}
\frac{d}{dt}(\ga , \z , \eta ) = (\eta, F(\ga,\eta,\z), G(\ga,\z))
\label{eq:ode}
\end{equation}
where 
\[ G(\ga,\z)=-(\mathrm{grad} \lop (\z \circ \gi)) \circ \ga ,
\]
\[ F(\ga,\z,\eta)= - ((\z \circ \gi)  \mathrm{div}(\eta \circ \gi)) \circ \ga) 
\]
with initial data $ (\ga_{0}, \eta_{0}, \z_{0})=(id_{x}, v_{0}(x), n_{0}(x)) $. 

In the proof of Theorem \ref{th1} we repeatedly use three standard results on Sobolev spaces: The Schauder ring property, Sobolev imbedding theorem (that we refer to as Sobolev lemma) and the composition lemma (see, for example, \cite{Stein} and \cite{Adams}).

\noindent {\bf Proof of Theorem \ref{th1}. }
If the map
\[ \begin{array}{rl}
\diff(\torm) \times \sobb(\torm,\mathbb{R}) \times \sob(\torm,\mathbb{R}^{m})
\rightarrow&\sob(\torm,\mathbb{R}^{m})
\times 
\sobb(\torm,\mathbb{R}) \times \sob(\torm,\mathbb{R}^{m}) \\ (\ga,\z,\eta) \mapsto&(\eta, F, G)
\end{array}
\]
is locally Lipschitz then by the fundamental theorem for ordinary
differential equations on Banach spaces \cite{Dieud} there is a unique solution 
\[ (\ga ,\eta ,\z) \in \diff(\torm) \times \sob(\torm,\mathbb{R}^{m}) \times \sobb(\torm,\mathbb{R}) 
\] 
to the problem (\ref{eq:ode}) for $s>m/2+1$ with initial data
\[ \z(0,x)=n_{0}(x) , \ \eta(0,x)=v_{0}(x) , \ \ga(0,x)=x.
\]

Note that the dependence of the solution of problem (\ref{eq:ode})
on initial data is smooth. However the map
$\ga\mapsto\gi$ on $\diff$ is continuous but not $C^{1}$. Therefore we only have continuous dependence on
initial data of the solution to the Cauchy problem for (mEP).   

By Proposition \ref{pro:ode_iap} the proof of Theorem \ref{th1} is reduced to showing that the maps
\begin{align}
\ga \mapsto &F(\ga,\eta,\z) \in L(\sob(\torm,\mathbb{R}^{m}), \sobb(\torm,\mathbb{R})), \nonumber \\
\z \mapsto &F(\ga,\eta,\z) \in L(\sobb(\torm,\mathbb{R}),\sobb(\torm,\mathbb{R})), \nonumber \\
\eta \mapsto &F(\ga,\eta,\z) \in L(\sob(\torm,\mathbb{R}^{m}),\sobb(\torm,\mathbb{R})), \label{eq:lala} \\ 
\ga \mapsto &G(\ga,\z) \in L(\sob(\torm,\mathbb{R}^{m}),\sob(\torm,\mathbb{R}^{m})), \nonumber \\
\z \mapsto &G(\ga,\z) \in L(\sobb(\torm,\mathbb{R}),\sob(\torm,\mathbb{R}^{m})) \nonumber
\end{align}
are locally Lipschitz in $\ga, \eta$ and $\z$ (uniformly with respect to the remaining variables). 

In the following estimates, the subscripts $\ga$ and $\z$ of a constant indicate the dependence of the constant on
$\| \ga \|_{\sob}$ and $\| \z \|_{\sobb}$ respectively. 
\\

{\bf $ \ga \mapsto G(\ga,\z) $ is locally Lipschitz:} 

Let $\ga_{1},\ga_{2}\in\diff(\torm)$ and $\bar{\z}\in\sobb(\torm,\reel)$. By the composition lemma,
\[ \| \mathrm{grad}\lop(\bar{\z}\circ\gi_{1})\circ\ga_{1}-\mathrm{grad}\lop(\bar{\z}\circ\gi_{2})\circ\ga_{2} \|_{\sob}
\]
\[ \leq C_{\ga_{2}} \| \mathrm{grad}\lop(\bar{\z}\circ\gi_{1})\circ\ga-\mathrm{grad}\lop(\bar{\z}\circ\gi_{2})\|_{\sob}
\]
where $\ga=\ga_{1} \circ\gi_{2}$. Let $\z=\bar{\z}\circ\ga_{2}^{-1}$. Then it is enough to show that the
following estimate holds:
\begin{equation} 
\| \mathrm{grad} \lop (\z \circ \gi )\circ \ga - \mathrm{grad}  \lop \z \|_{\sob} \leq C_{\ga,\z} \| \ga - id_{x} \|_{\sob}. 
\label{eq:cross}
\end{equation}

Next we show that (\ref{eq:cross}) holds for $s>m/2+1$. The left side of the inequality (\ref{eq:cross}) has the following form
\begin{align}
& \| \grad \lop (\z \circ \gi )\circ \ga - \grad \lop \z\|_{\sob} \nonumber \\ 
& \simeq \| \grad\lop (\z \circ \gi )\circ \ga - \grad \lop \z\|_{L^{2}} \label{le2} \\
& + \| \div (\grad\lop (\z \circ \gi )\circ \ga - \grad \lop \z)\|_{\sobb}. \label{es-1}
\end{align}

We first estimate the $L^{2}$ term in (\ref{le2}). For any $r>m/2+\sigma$ we have the Sobolev
imbedding  into a H\"{o}lder space
\[ H^{r}(\torm) \hookrightarrow C^{\sigma}(\torm)
\]
with a bound
\begin{equation}
|u(x)-u(y)|\leq C \| u\|_{H^{r}} |x-y|^{\sigma},
\label{eq:yardim}
\end{equation}
for any $x,y \in \tor$. Therefore using the composition lemma and applying (\ref{eq:yardim}) with $r=s-1$ 
and $\sigma = (s-1-m/2)/2$ we find
\begin{align*}
& \| \grad\lop (\z \circ \gi)\circ \ga - \grad\lop\z\|_{L^{2}} \\
& \leq \| \grad\lop (\z \circ \gi)\circ \ga - \grad\lop(\z\circ \gi)\|_{L^{2}} + \| \grad\lop (\z \circ \gi) -
\grad\lop\z\|_{L^{2}} \\
& \leq C_{\ga} \| \z\|_{\sobb} \left( \| \ga - id_{x} \|_{\sob}^{\sigma +1} + \| \gi - id_{x} \|_{\sob}^{\sigma} \right).
\end{align*} 
Adding and subtracting the appropriate terms we estimate the $\sobb$ term in (\ref{es-1}) by the sum
\begin{align}
\leq & \left\| \sum_{i=1}^{n}\sum_{j=1}^{n}\dj\di\lop (\z \circ\gi)\circ\ga (\di\ga_{j}-\delta_{i}^{j})\right\|_{\sobb} \label{ro1}\\
 & + \| \del\lop(\z\circ\gi)\circ\ga- \del\lop\z \|_{\sobb}. \label{ro2}
\end{align}
Using Schauder ring property and composition lemma the first summand in (\ref{ro1}) is bounded by
\begin{align*}
&\leq \| \z\|_{\sobb} \| \ga -id_{x}\|_{\sob}.
\end{align*}
In order to estimate the second summand (\ref{ro2}) we add and subtract the terms $\lop (\z\circ\gi)\circ\ga$ and $\lop\z$. After cancellations we obtain
%\begin{align*}
%& = \| \lop (\z\circ\gi)\circ\ga - \lop\z \|_{\sobb}.
%\end{align*}
%By the triangle inequality we have
\begin{align}
& =\| \lop (\z\circ\gi) \circ\ga - \lop \z \|_{\sobb} \nonumber \\
&\leq \| \lop (\z\circ\gi) \circ\ga - \lop
(\z\circ\gi)\|_{\sobb}  \label{eq:ilki}\\
& + \| \lop (\z\circ\gi) - \lop \z \|_{\sobb}. \label{eq:ikinc}
\end{align}

Let $u$ be 
$\lop(\z\circ\gi)$. Then we have
\begin{align*}
\| \lop (\z\circ\gi) \circ\ga - \lop (\z\circ\gi)\|_{\sobb}& = \left\| \int_{x}^{\ga(x)}  
D u(y)dy \right\|_{\sobb} \\
& \leq \| u\|_{C^{1}} \|\ga -id_{x}\|_{\sobb}.
\end{align*}
Using the Sobolev lemma with the composition lemma we obtain the estimate 
\begin{align}
& \leq \|  u \|_{\sob} \|\ga - id_{x}\|_{\sobb} \nonumber \\
& \leq C_{\ga} \| \z\|_{\sobb} \|\ga -id_{x}\|_{\sob} \label{klik}
\end{align}
for (\ref{eq:ilki}).

For $s > m/2+2$ the term (\ref{eq:ikinc}) can easily be estimated like (\ref{eq:ilki}). 
For $m/2+1<s\leq m/2+2$ we
first observe that the estimate
\begin{equation} \| \lop (\z\circ\gi) - \lop \z \|_{\sobb} \leq \| \z\circ\gi -  \z \|_{L^{2}}
\label{leyla1}
\end{equation}
holds for $s-3\leq 0$ and then using (\ref{eq:yardim}) as before, we obtain the following estimate for
(\ref{eq:ikinc})
\begin{align}
& \leq C_{\ga} \| \z\|_{\sobb} \| \gi -id_{x}\|_{\sob}^{\sigma} \label{leyla2}
\end{align}
where $\sigma$ is equal to $(s-1-m/2)/2>0$.
However the assumption $s-3\leq 0$ does not follow from $s\leq m/2+2$ if $m\geq 3$. Nevertheless one can use the following inductive argument until $s-(2k+1)\leq 0$ (it ends in finitely many steps since $s\leq m/2+2$). 

If $s-3>0$ we split (\ref{eq:ikinc}) as in (\ref{le2})-(\ref{es-1}). The $L^{2}$ part can be estimated as for $s-3\leq 0$. The $H^{s-2}$ part 
\begin{align}
& \| D(\lop(\z\circ\gi)\circ\ga-\lop\z) \|_{H^{s-2}} \nonumber \\
& = \|\grad\lop(\z\circ\gi)\circ\ga . D\ga -\grad\lop\z \|_{H^{s-2}} \label{eq:wasa}
\end{align}
is bounded by
\begin{align}
& \leq \|\grad\lop(\z\circ\gi)\circ\ga . (D\ga -1) \|_{H^{s-2}} \label{eq:bdokuz} \\
& \ +\|\grad\lop(\z\circ\gi)\circ\ga -\grad\lop\z \|_{H^{s-2}} \label{eq:dokuz}
\end{align}
Here the first summand (\ref{eq:bdokuz}) is estimated using the Schauder ring property with the composition lemma
\[ \|\grad\lop(\z\circ\gi)\circ\ga . (D\ga -1) \|_{H^{s-2}} \leq C_{\ga}\|\z\|_{\sobb} \|\ga-id_{x}\|_{\sob}.
\] 
For the second summand (\ref{eq:dokuz}) we use the steps (\ref{eq:cross})-(\ref{klik}) to reduce it to estimating
\begin{equation} 
\|\lop (\z\circ\gi)-\lop\z \|_{H^{s-3}}. 
\label{pom}
\end{equation}
If $s-5 \leq 0$ we proceed as in (\ref{leyla1})-(\ref{leyla2}). Otherwise we repeat the steps (\ref{eq:wasa})-(\ref{pom}).

{\bf $ \ga \mapsto F(\ga,\z,\eta) $ is locally Lipschitz:}

Let $\ga_{1},\ga_{2}\in\sob(\torm,\mathbb{R}^{m})$ and $\z\in\sobb(\torm,\reel), \eta\in\sob(\torm,\reel^{m})$.
Then by Schauder ring property we have 
\begin{align*}
& \| \z (\mathrm{div}(\eta\circ\gi_{1})\circ\ga_{1})-\z (\mathrm{div}(\eta\circ\gi_{2})\circ\ga_{2})\|_{H^{s-1}} \\
%& \leq\| \z \|_{H^{s-1}} \| \sum_{j} \dj(\eta_{j}\circ\gi_{1})\circ\ga_{1}-\dj(\eta_{j}\circ\gi_{2})\circ\ga_{2}\|_{H^{s-1}} \\
& \leq \|\z \|_{H^{s-1}} \| \sum_{j,m} \dm \eta_{j}\dj(\gi_{1})_{m}\circ\ga_{1}-\dm\eta_{j}
\dj(\gi_{2})_{m}\circ\ga_{2}\|_{H^{s-1}} \\
& \leq C_{\ga_{1},\ga_{2},\z} \sum_{j,m}\| \dm\eta_{j} \|_{H^{s-1}} \|
\dj(\ga_{1}^{-1})_{m}-\dj(\ga_{2}^{-1})_{m}\|_{H^{s-1}}.
\end{align*}
Using the Schauder ring property one more time we bound this term by
\begin{align*}
& \leq C_{\ga_{1},\ga_{2},\z,\eta} \| \eta\|_{\sob}\| \ga_{1}-\ga_{2} \|_{\sob}.
\end{align*}
and therefore $ \ga \mapsto F(\ga,\z,\eta) $ is locally Lipschitz.

It is straightforward to show that the second, third and fifth maps in (\ref{eq:lala}) are uniformly Lipschitz using properties of Sobolev spaces.
This completes the proof of Theorem \ref{th1}.

\hfill $ \Box $

Next we observe that the Cauchy problems for the equation (mEP) and the Euler equations of an 
incompressible fluid are not only similar for low regularity (Sobolev class) data but also for high regularity (analytic) data.
 
\section{Analytic regularity} 

In this section we give a proof of theorem \ref{th:iap_analy_reg} that states the analytic regularity (i.e., 
existence and uniqueness of analytic solutions for analytic initial data) of the Cauchy problem for (mEP).

Our approach is motivated by the work of M.S. Baouendi and C. Goulaouic \cite{BG} who studied analytic regularity of the
Cauchy problem for Euler equations of incompressible fluids.

The proof of theorem \ref{th:iap_analy_reg} relies on a contraction argument in 
a decreasing scale
of Banach spaces $X_{s}$ (i.e. if $s'<s$ implies $X_{s}\subset X_{s'}$ and $|||\cdot |||_{s'}\leq |||\cdot |||_{s}$).
 
For $s>0$, let the spaces $E_{s}$ be defined as
\[ E_{s}=\left\{ u \in C^{\infty}(\torm) : \int_{\torm} u \ dx=0 \ and \ ||| u |||_{s}=\sup_{|k|\geq 0} \frac{\| \dx^{k}
u \|_{H^{\sigma}}s^{|k|}}{k! / (|k|+1)^{2}} < \infty \right\},
\]
where $\sigma$ is any integer such that $\sigma >1+m/2$ and let $X_{s}$ be given by the Cartesian product $E_{s}\times E_{s}$. The norm $|||\cdot |||_{X_{s}}$ can be chosen to be any of the
standard product norms on $E_{s}\times E_{s}$. The following lemma states the ring property for the spaces $E_{s}$.   
\begin{lema}
Let $0<s<1$. There is a constant $c>0$ which is independent of $s$ such that we have
\[ ||| uv |||_{s}\leq c |||u|||_{s} |||v|||_{s}
\]
for any $u,v \in E_{s}$.
\label{lema:product}
\end{lema}

We omit the proof of lemma \ref{lema:product} (see \cite{HM1} for the case $m=1$). 
First we rewrite the equation (mEP) in a more convenient form.
Let $n$ and $v$ be denoted by $u_{1}$ and $u_{2}$ respectively. Then we can write the equation
(mEP) in terms of $(u_{1},u_{2})$ as 
\begin{equation} \begin{array}{l}
\dt u_{1}=F_{1}(u_{1},u_{2}):= \P_{2}(u_{1}u_{2}), \\ 
\dt u_{2}=F_{2}(u_{1},u_{2}):= \P_{4}(u_{2})u_{2}+\P_{1}\P_{3}u_{1}
\end{array}
\label{eq:syst_iap}
\end{equation}
where 
\[\begin{array}{ll}
 \P_{1}(n):=-\mathrm{grad}(n), & \P_{3}(n):=\lop n, \\
 \P_{2}(v):=-\mathrm{div}(v), & \P_{4}(u)v:=-(\gr_{v}u)=-(Du)v.
\end{array}
\]

The following lemmas give the suitable bounds on these operators to prove Theorem
\ref{th:iap_analy_reg}. 

\begin{lema}
For $0<s'<s<1$, we have   
\[ ||| \P_{1}n|||_{s'} \leq \frac{C}{s-s'} ||| n |||_{s}, \]
\[ |||\P_{2}v|||_{s'}\leq \frac{c}{s-s'}|||v|||_{s}.
\]
\label{lema:Pone}
\end{lema}
\pf By the definition of $|||.|||_{s}$, we have
\[ ||| \P_{1}n|||_{s'} = \sup_{|k|\geq 0} \frac{\| \partial^{k} \P_{1}n \|_{H^{\sigma}} {s'}^{|k|}}{k!/(|k|+1)^{2}}.
\]
The $H^{\sigma}$ norm on the right hand side can be written in the local coordinates up to a constant as
\[ \| \partial^{k}(\mathrm{grad}(n))\|_{H^{\sigma}(\torm, \reel^{m})} \simeq \sum_{j=1}^{m} \|
\partial^{k}\partial_{j}n\|_{H^{\sigma}(\torm,\reel)}.
\]
Then we have the estimates
\begin{align*}
||| \P_{1}n|||_{s'}& \leq c \sup_{|k|\geq 0}\sup_{|\beta|=1} \frac{\| \partial^{k+\beta} n \|_{H^{\sigma}}
{s'}^{|k|}}{k!/(|k|+1)^{2}} \\
& = c \sup_{|k|\geq 0}\sup_{|\beta|=1} \frac{\| \partial^{k+\beta} n \|_{H^{\sigma}}
{s}^{|k|+1}}{(k+\beta)!/(|k|+2)^{2}}\frac{(k+\beta)!}{(|k|+2)^{2}}\frac{s'^{|k|}}{s^{|k|+1}}\frac{(|k|+1)^{2}}{k!}.
\end{align*} 
Note that 
\[ \sup_{|k|\geq 0}\sup_{|\beta|=1} \frac{\| \partial^{k+\beta} n \|_{H^{\sigma}} {s}^{|k|+1}}{(k+\beta)!/(|k|+2)^{2}} \leq
||| n |||_{s}.
\]
Therefore we have
\begin{equation}
||| \P_{1}n|||_{s'}\leq m ||| n |||_{s}\sup_{|k|\geq 0}\sup_{|\beta|=1}
\frac{s'^{|k|}}{s^{|k|+1}}\frac{(k+\beta)!}{k!}\left( \frac{|k|+1}{|k|+2}\right)^{2}.
\label{eq:lale}
\end{equation}
Note also that 
\[ \sup_{|\beta|=1}\frac{(k+\beta)!}{k!}= \sup_{1\leq i \leq m} (k_{i}+1) \leq |k|+1.
\]
Then it follows from the formula (\ref{eq:alg}) that the inequality
\begin{equation} 
\sup_{|k|\geq 0}\sup_{|\beta|=1}
\frac{s'^{|k|}}{s^{|k|+1}}\frac{(k+\beta)!}{k!}\left( \frac{|k|+1}{|k|+2}\right)^{2} \leq \frac{C}{s-s'}
\label{eq:gul}
\end{equation}
holds. 
By (\ref{eq:lale}) and (\ref{eq:gul}) we obtain
\[ |||\P_{1}n|||_{s'}\leq \frac{C}{s-s'}|||n|||_{s}.
\]
The estimate for $\P_{2}$ follows similarly.

\hfill $ \Box $

\begin{lema}
For any $0<s<1$, the estimate 
\[ ||| \P_{3}(u)|||_{s} \leq ||| u |||_{s} \]
holds if $u \in E_{s}$.
\label{lema:P3}
\end{lema}

\begin{lema}
For $0<s'<s<1$, we have 
\[ |||\P_{4}(u)v|||_{s'}\leq \frac{c}{s-s'} |||v|||_{s'} |||u|||_{s}.
\]
\label{lema:P4}
\end{lema}
\pf We write $\P_{4}(u)v$ in terms of the linear operator $Du$ as $\P_{4}(u)v=(Du)v$. Then by Lemma \ref{lema:product} we
have
\begin{align*}
|||\P_{4}(u)v|||_{s'} & =|||(Du)v|||_{s'} \\
& = |||\sum_{i,j=1}^{m} v_{i}\partial_{i}u_{j}|||_{s'} \\
& \leq \sum_{i,j=1}^{m} |||v_{i}|||_{s'} |||\partial_{i} u_{j}|||_{s'}.
\end{align*}
Now we reduce the proof to the case that we handled in the proof of lemma \ref{lema:Pone}:
\begin{align*}
|||\P_{4}(u)v|||_{s'} & \leq c |||v|||_{s'} \sum_{i,j=1}^{m} |||\partial_{i}u_{j}|||_{s'} \\
& \leq c |||v|||_{s'} \sum_{j=1}^{m} \sup_{|k|\geq 0}\sup_{|\beta|=1}\frac{\| \partial^{k+\beta} u_{j}\|_{H^{\sigma}}
s'^{|k|}}{k!/(|k|+1)^{2}} \\
& \leq c |||v|||_{s'} |||u|||_{s} \sup_{|k|\geq 0}\sup_{|\beta|=1}\frac{s'^{|k|}}{s^{|k|+1}}\frac{(k+\beta)!}{k!}\left(
\frac{|k|+1}{|k|+2}\right)^{2}.  
\end{align*}
Clearly, to finish the proof, it is enough to show that
\begin{equation} 
\frac{{s'}^{k}}{s^{k+1}}\left( \frac{k+1}{k+2} \right)^{2} (k+1) \leq \frac{1}{s-s'}. 
\label{eq:alg}
\end{equation}
Let $s'=\lambda s$, $0< \lambda <1$ and $f(\lambda)=(k+1)(1-\lambda)\lambda^{k}$. Then,
\begin{align*}
\frac{{s'}^{k}}{s^{k+1}}\left( \frac{k+1}{k+2} \right)^{2} (k+1) = & \frac{\lambda^{k} s^{k}}{s^{k+1}}\left(
\frac{k+1}{k+2} \right)^{2} (k+1) \\
=& \frac{1}{s(1-\lambda)} f(\lambda) \left(
\frac{k+1}{k+2} \right)^{2}.
\end{align*}

For $k=0$ it is clear that $f(\lambda)\leq 1$. For $k\geq 1$ the function  $f(\lambda)=(k+1)(1-\lambda)\lambda^{k}$ is
continuous in the interval $0<\lambda <1$ and it has zeros at the endpoints of the interval $[0,1]$ and a maximum at
$\lambda=\frac{k}{k+1}$ such that
$f(\frac{k}{k+1}) =\left(\frac{k}{k+1}\right)^{k}<1$. Then we have  
\[ (\lambda^{k}-\lambda^{k+1}) \left( \frac{k+1}{k+2} \right)^{2} (k+1) \leq 1 .
\]
Therefore
\[ \frac{\lambda^{k} s^{k}}{s^{k+1}}\left( \frac{k+1}{k+2} \right)^{2} (k+1) \leq \frac{1}{s(1-\lambda)}
\]
and the formula (\ref{eq:alg}) holds.
Using the formula (\ref{eq:alg}) we obtain the desired estimate
\[ |||\P_{4}(u)v|||_{s'} \leq \frac{c}{s-s'} |||v|||_{s'} |||u|||_{s}.
\] 

\hfill $ \Box $

Now we are ready to prove theorem \ref{th:iap_analy_reg}.

{\bf proof of theorem \ref{th:iap_analy_reg}.} We refer to the version of the abstract Cauchy-Kowalevski
theorem in \cite{Nis}. We only need to verify
the first two conditions of this theorem since the map $F(u_{1},u_{2})$ does not depend on $t$
explicitly.

Clearly, $t\longmapsto F(t,u(t))=(F_{1}(u_{1},u_{2}),F_{2}(u_{1},u_{2}))$ is holomorphic if $t\longmapsto u_{1}(t)$ and 
$t\longmapsto u_{2}(t)$ are both holomorphic. We only need to show that $F_{1}(u_{1},u_{2})$ and $F_{2}(u_{1},u_{2})$ are in
$ E_{s'}$ if $u_{1},u_{2}\in E_{s}$. 
We begin with estimates on $F_{1}$.

By Lemma \ref{lema:Pone} and Lemma \ref{lema:product}, we have
\[ |||F_{1}(u_{1},u_{2})|||_{s'}  = ||| \P_{2}(u_{1}u_{2}) |||_{s'} \leq \frac{c}{s-s'} |||u_{1}|||_{s} |||u_{2}|||_{s}.
\]

Similarly, for $F_{2}$, using the lemmas \ref{lema:Pone}, \ref{lema:P3} and \ref{lema:P4} we have
\[ |||F_{2}(u_{1},u_{2}) |||_{s'}  = |||\P_{4}(u_{2})u_{2} + \P_{1}\P_{3}u_{1} |||_{s'}\leq \frac{c}{s-s'} \left( |||u_{2}|||_{s}^{2} + |||u_{1}|||_{s} \right).
\]

We proceed to establish the second condition of the abstract Cauchy-Kowalevski theorem. We will show that for some
$c$ independent of $t$,
\[  ||| F_{1}(u_{1},u_{2})- F_{1}(v_{1},v_{2})|||_{s'} \leq \frac{c}{s-s'} |||u-v|||_{X_{s}}
\]
and similarly
\[  ||| F_{2}(u_{1},u_{2})- F_{2}(v_{1},v_{2})|||_{s'} \leq \frac{c}{s-s'} |||u-v|||_{X_{s}}
\]
where $u=(u_{1},u_{2})$ and $v=(v_{1},v_{2})$.

To obtain the first estimate above, after applying Lemma \ref{lema:Pone}, we add and subtract the term $u_{1}v_{2}$ and use
Lemma \ref{lema:product}: 
\begin{align*}
||| F_{1}(u_{1},u_{2})-F_{1}(v_{1},v_{2}) |||_{s'} & =|||\P_{2}(u_{1}u_{2}-v_{1}v_{2}) |||_{s'} \\
%& \leq \frac{c}{s-s'} ||| u_{1}u_{2}-v_{1}v_{2} |||_{s} \\
%& \leq \frac{c}{s-s'} \left( |||u_{1}(u_{2}-v_{2})|||_{s} + |||(u_{1}-v_{1})v_{2}|||_{s} \right) \\
& \leq \frac{c}{s-s'} \left( |||u_{1}|||_{s}|||u_{2}-v_{2}|||_{s} + |||u_{1}-v_{1}|||_{s}|||v_{2}|||_{s} \right).
\end{align*}
Then, assuming that $|||u|||_{s}<R$ and $|||v|||_{s}<R$, we have
\[ ||| F_{1}(u_{1},u_{2})- F_{1}(v_{1},v_{2})|||_{s'}  \leq \frac{c}{s-s'} |||u-v|||_{X_{s}}.
\]
To estimate the $F_{2}$ component, we use lemmas \ref{lema:Pone} and \ref{lema:P3}:
\begin{align}
|||F_{2}(u_{1},u_{2})-F_{2}(v_{1},v_{2})|||_{s'} =
|||(\P_{4}u_{2})u_{2}-(\P_{4}v_{2})v_{2}+\P_{1}\P_{3}(u_{1}-v_{1})|||_{s'} \nonumber \\
\leq |||(\P_{4}u_{2})u_{2}-(\P_{4}v_{2})v_{2}|||_{s'} + \frac{c}{s-s'}||| u_{1}-v_{1} |||_{s}.
\label{eq:pfour}
\end{align}
Note that
\[ \P_{4}(u_{2})u_{2}-\P_{4}(v_{2})v_{2}=\P_{4}(u_{2})(u_{2}-v_{2})+\P_{4}(u_{2}-v_{2})v_{2}.
\]
Using Lemma \ref{lema:P4} and the above identity, (\ref{eq:pfour}) implies
\begin{align*}
|||F_{2}(u_{1},u_{2})-F_{2}(v_{1},v_{2})|||_{s'} \leq \frac{c}{s-s'}& \left( |||u_{2}|||_{s}
|||u_{2}-v_{2}|||_{s'} + \right. \\
& \left. + |||u_{2}-v_{2}|||_{s}|||v_{2}|||_{s'} + |||u_{1}-v_{1}|||_{s} \right).
\end{align*}
Therefore the estimate
\[ |||F_{2}(u_{1},v_{1})-F_{2}(u_{2},v_{2}) |||_{s'} \leq \frac{c}{s-s'}|||u-v|||_{X_{s}}
\]
holds. This completes the proof of theorem \ref{th:iap_analy_reg}.

\hfill $\Box$

\section{Bihamiltonian structure and integrability}

A number of partial differential equations that describe fluid motion can be derived as equations for
geodesics on various infinite dimensional Lie groups. For instance, the Euler equation for ideal
incompressible fluid flow is the geodesic equation on the group of volume-preserving diffeomorphisms of
a Riemannian manifold $M$ with a right invariant metric given by the $L^{2}$ inner product on the
tangent space at the identity of the group \cite{EMa}. Other examples are 
\begin{itemize}
\item Korteweg-de Vries equation and
Camassa-Holm equation on the Bott-Virasoro group (see for example \cite{OK} and \cite{Mis3}), 
\item ideal incompressible MHD (magnetohydrodynamics) on the semidirect product of volume preserving
diffeomorphisms with the divergence free vector fields, 
\item Hunter-Saxton equation on the homogeneous space of all diffeomorphisms of the unit circle
modulo the rotations \cite{KhMi}, etc.
\end{itemize}

In contrast with all the examples we gave above, the energy of the modified Euler-Poisson equation (mEP) is not a
quadratic form, therefore it can not be interpreted as a Riemannian metric. However, there still is a
variational problem on the cotangent space of the configuration space of this equation. Here we derive 
the equation (mEP) from this variational problem. Note that all the computations
that follow are formal.

Let
$\mathfrak{g}$ be a Lie algebra with the bracket operation $[\cdot ,\cdot]$ and $\mathfrak{g}^{*}$ be its dual given by
the pairing
\[ \langle\cdot ,\cdot\rangle :\mathfrak{g}^{*}\times\mathfrak{g} \rightarrow \reel .
\] 
Then $\mathfrak{g}^{*}$ with the Lie-Poisson bracket defined by 
\[ \{ F,G\}(m)=\langle m , \left[ \frac{\delta F}{\delta m}, \frac{\delta G}{\delta m}\right] \rangle
\]
for $F,G:\mathfrak{g}^{*}\rightarrow \reel$ is a Poisson manifold.

Hamilton's equations for $H\in C^{\infty}(P)$ on a Poisson manifold $(P, \{\cdot ,\cdot\})$ are given by
\begin{equation}
\dot{F}=\{ F,H\} \ \mathrm{for} \ \mathrm{all} \ F\in C^{\infty}.
\label{eq:hamilton1}
\end{equation}
Given a Lie-Poisson manifold, Hamilton's equations (\ref{eq:hamilton1}) can be written as 
\begin{align*}
\dot{F}(m)&=\{ F,H\}(m) \\
\langle \frac{\delta F}{\delta m}, \dt m \rangle &= -\langle m, \left[ \frac{\delta H}{\delta m},\frac{\delta F}{\delta m}\right] \rangle \\
& =-\langle ad^{*}_{\delta H/\delta m}m, \frac{\delta F}{\delta m} \rangle. 
\end{align*}
Therefore 
\begin{equation}
\dt m = - ad^{*}_{\delta H/\delta m}m
\label{eq:hamilton2}
\end{equation}
is an equivalent formulation of the Hamilton's equation on a Lie-Poisson manifold. 

Here we exploit the tools and techniques used to study the Hamiltonian formulation of the
Euler equations for a compressible fluid \cite{MRW} to show that the modified Euler-Poisson equation (mEP) can be derived as a
Hamiltonian equation.  

On the Cartesian product space $\mathrm{Diff}(\torm) \times C^{\infty}(\torm)$ of the group of diffeomorphisms of
$\torm$ and the vector space $C^{\infty}(\torm)$ of all smooth functions on $\torm$, the operation
\[ (\phi,a)\circ(\psi,b)=(\phi\circ\psi, a\circ\psi^{-1}+b),
\]
called the semidirect product induces a Lie group structure. We denote this group by 
\[ G=\mathrm{Diff}(\torm)\ltimes C^{\infty}(\torm)
\]
following the conventional notation for semidirect product spaces.
The corresponding Lie algebra is the space $\mathfrak{g}=\mathrm{Vect}(\torm)\ltimes C^{\infty}(\torm)$ with the
bracket
\[ [(v,a),(w,b)]=([v,w], \mathcal{L}_{v} b- \mathcal{L}_{w} a)
\] 
where $v,w\in \mathrm{Vect}(\torm)$ and $a,b\in C^{\infty}(\torm)$. Here $[v,w]$ is the usual commutator of vector fields
on $\torm$ and $\mathcal{L}_{w}a$ is the Lie derivative of $a$ in the direction of $w$ and is given  by $ \mathcal{L}_{w}
a=\left. \frac{d}{ds}\right|_{s=0}(a\circ \z^{s})$ where $\z^{s}$ is any curve on $\mathrm{Diff}(\torm)$ such that
$\z^{s}|_{s=0}=id$ and $\left. \frac{d}{ds}\right|_{s=0}\z^{s}=w $. Note that in this setting the composition of the
diffeomorphisms is the group operation on $\mathrm{Diff}(\torm)$, the composition of a smooth function with a
diffeomorphism $a\circ \gi$ is the natural action of the diffeomorphism $\ga$ on the function $a$. In general, the
semidirect product structure on the Cartesian product of a Lie group and a vector space on which the group acts is defined
using the group operation and the action of the group on the vector space (see \cite{ArKh}, \cite{MRW}). 

In this context, a Hamiltonian formulation of the modified Euler-Poisson equation (mEP) can be stated as follows:

\begin{teor}
The modified Euler-Poisson equation (mEP) is a Hamiltonian equation on $\mathfrak{g}^{*}$ with respect to the
linear Lie-Poisson structure and the energy function
\begin{equation}
H(M,n)=-\int_{\torm}\left( \frac{1}{2n}\langle M,M\rangle +\Phi(n) \right) dx
\label{eq:Hmulti}
\end{equation}
where $M=nv\in \mathrm{Vect}(\torm)$ and $\Phi'(n)=\lop(n)$.
\label{th:hamilt}
\end{teor}

\pf We want to derive the equations for $v$ and $n$ from
\begin{equation} 
\dt m=-ad^{*}_{\delta H/\delta m}m
\label{eq:ad}
\end{equation}
where $m=(M,n)=(nv,n)$ and $(v,n) \in \mathfrak{g}^{*}$.

The variational derivative $\frac{\delta H}{\delta m}$ is given by $\left( \frac{\delta
H}{\delta M}, \frac{\delta H_{1}}{\delta n}\right)$ with
\begin{equation}
\frac{\delta H_{1}}{\delta M}= -M/n=-v,
\label{eq:h1}
\end{equation}
\begin{equation}
\frac{\delta H_{1}}{\delta n}= \frac{1}{2n^{2}} \langle M,M \rangle -\Phi'(n)= \frac{1}{2}\langle v,v \rangle -\Phi'(n).
\label{eq:h2}
\end{equation}

Evaluating equation (\ref{eq:ad}) on an arbitrary pair $(w,b)\in \mathfrak{g}$, we obtain
\[ \langle \dt m, (w,b)\rangle =-\langle ad^{*}_{\delta H_{1}/\delta m}m, (w,b)\rangle .
\]
Then by the definition of the coadjoint operator $ad^{*}$ and the bracket on $\mathfrak{g}$, we have 
\begin{align*}
\langle \dt m, (w,b)\rangle &=-\langle (M,n),\left[ \left( \frac{\delta H_{1}}{\delta M}, \frac{\delta
H_{1}}{\delta n}\right), (w,b) \right] \rangle \\   
& = \langle (M,n),\left( \left[ w, \frac{\delta H_{1}}{\delta M} \right], \mathcal{L}_{w}\frac{\delta
H_{1}}{\delta n}-\mathcal{L}_{\frac{\delta H_{1}}{\delta M}}b \right)\rangle.
\end{align*}
In what follows we identify the dual space $\mathfrak{g}^{*}$ with the algebra $\mathfrak{g}$ using the pairing
$\langle \cdot , \cdot \rangle $ on $\mathfrak{g}^{*}\times \mathfrak{g}$ given by
\begin{equation} 
\langle (v,a),(w,b)\rangle =\int_{\torm} \langle v,w \rangle \ dx + \int_{\torm} ab \ dx .
\label{eq:couple}
\end{equation} 
Then substituting $nv$ for $M$ and using (\ref{eq:h1}) and (\ref{eq:h2}), we obtain 
\begin{align}
\langle \dt m, (w,b)\rangle = & \int_{\torm} \langle [w,-v], nv \rangle dx + \nonumber \\
 & + \int_{\torm} n\left(
\mathcal{L}_{w}(\frac{1}{2}\langle v,v \rangle-\Phi'(n))-\mathcal{L}_{-v}b \right) dx. \nonumber
\label{eq:!}
\end{align}
By the definition of the bracket $[\cdot , \cdot ]$ on $\mathrm{Vect(\torm)}$, we have 
\begin{align*}
\langle \dt m,(w,b)\rangle = & \int_{\torm}\langle (v.\gr)w,nv \rangle - \langle (w.\gr)v,nv \rangle dx \\ 
 & + \int_{\torm} n\left(
\mathcal{L}_{w}(\frac{1}{2}\langle v,v \rangle-\Phi'(n))-\mathcal{L}_{-v}b \right) dx.
\end{align*}
Furthermore, we can compute the Lie derivatives on the right hand side and write the above equality as follows
\begin{align*} 
\langle \dt m,(w,b)\rangle = & \int_{\torm}\langle (v.\gr)w,nv \rangle - \langle (w.\gr)v,nv \rangle dx +\\ 
 & + \int_{\torm} \langle nw, \mathrm{grad}(\frac{1}{2}\langle v, v \rangle-\Phi'(n) ) \rangle + \langle nv,
\mathrm{grad}(b) \rangle dx.
\end{align*}
Using the identities
\[ \langle v, \mathrm{grad} \langle v, w \rangle \rangle=\langle w, (v.\gr)v \rangle + \langle v, (v.\gr)w \rangle
\]
and
\[ \langle w, \frac{1}{2}\mathrm{grad}\langle v,v \rangle \rangle =\langle v, (w.\gr)v \rangle
\]
we obtain
\begin{align*}
\langle \dt m,(w,b)\rangle =&\int_{\torm}\langle\mathrm{grad}\langle v,w \rangle , nv \rangle - \langle
(v.\gr)v, nw \rangle - \\ 
& \ \ - \langle nw, \mathrm{grad}(\Phi'(n) ) \rangle + \langle nv, \mathrm{grad}(b) \rangle dx.
\end{align*}
Integrating by parts the first and the last summands on the right hand side, we obtain
\[ \langle \dt m,(w,b)\rangle = \int_{\torm} - \mathrm{div}(nv) \langle v,w \rangle - \langle (v.\gr)v, nw
\rangle - \langle nw, \mathrm{grad}(\Phi'(n) ) \rangle -b \mathrm{div}(nv) dx.
\]
Note that 
\[ \langle \dt m,(w,b)\rangle = \langle ((\dt n)v+n(\dt v), \dt n), (w,b) \rangle .
\]
Then by (\ref{eq:couple}), we have
\[ \dt n=-\mathrm{div}(nv),
\]
\[ (\dt n)v+n(\dt v)=-\mathrm{div}(nv)v - n(v.\gr)v - n \mathrm{grad}(\Phi'(n))
\]
which is equivalent to
\begin{align*} 
& \dt n=-\mathrm{div}(nv), \\
& \dt v=-(v.\gr)v-\mathrm{grad}(\Phi'(n)).
\end{align*}

\hfill $\Box$
\\

Note that for one space dimension ($m=1$) the hamiltonian $H_{1}$ in (\ref{eq:Hmulti}) is given by 
\begin{equation}
H_{1}=\int\frac{1}{2}\left(v^{2}n+(\lop\dx n)^{2}+(\lop n)^{2}\right)dx
\end{equation}
in terms of $v$ and $n$. Then using the differential operator $\D_{1}$ that is defined as
\[ \D_{1}=\left( \begin{array}{rr}
0 & -\partial \\
-\partial & 0
\end{array}\right),
\]
one can rewrite the equation (mEP) in the hamiltonian form
\begin{equation}
\dt\left( \begin{array}{c}
v \\
n
\end{array}
\right) = \D_{1}\left(\begin{array}{c}
\delta H_{1}/\delta v \\
\delta H_{1}/\delta n
\end{array} \right).
\label{eq:poi1}
\end{equation}
The Poisson bracket induced by the matrix differential operator $\D_{1}$ is given by 
\[ \{ F,H \}(v,n)=\langle (v,n),\left( \left[\left( \fv ,\fn \right) ,\left(\hv ,\hn \right)\right]\right)\rangle \]
\[ = \int \left(\fn \partial_{x}\hv -\hn \partial_{x}\fv \right) dx.
\]
Another conserved quantity for (mEP) is 
\begin{equation}
H_{2}=\int nv \ dx.
\end{equation}
For one space dimension ($m=1$) we can use $H_{2}$ to write (mEP) in yet another form as 
\begin{equation}
\dt\left( \begin{array}{c}
v \\
n
\end{array}
\right) = \D_{2}\left(\begin{array}{c}
\delta H_{2}/\delta v \\
\delta H_{2}/\delta n
\end{array} \right)
\label{eq:poi2}
\end{equation}
where $\D_{2}$ is defined as
\[ \D_{2}=\left( \begin{array}{cc}
\lop\dx & -\dx v \\
\dx v & -(n\partial+\partial n)
\end{array}\right).
\]
We prove theorem \ref{th:integrable} by showing that (\ref{eq:poi1}) and (\ref{eq:poi2}) are Hamiltonian forms of the modified Euler-Poisson equation (mEP) and that the induced Poisson structures are compatible. 

{\bf proof of theorem \ref{th:integrable}.}
The matrix differential operator $\D_{1}$ is skew-adjoint and does not depend on $v$ nor $n$ nor any of their derivatives, therefore the bracket given by $\D_{1}$ satisfies the Jacobi identity hence is indeed a Poisson bracket. 
%(as we described earlier, it is the Poisson bracket induced by the Lie bracket on the semidirect product Lie algebra %$Vect(\tor)\ltimes C^{\infty}$).

We can easily check that $\D_{2}$ is skew-adjoint as well:
\[ \int (\po , \pt)\D_{2}(\to , \tt)dx 
\]
\[= \int \po (\lop\dx\to - \tt \dx v)+ \pt (\to\dx v - n\dx \tt -\dx(n\tt))dx
\]
\[= -\int (\to\lop\dx\po - \pt\to\dx v+\po\tt\dx v-\tt\dx(\pt n)-n\tt\dx\pt) dx
\]
\[= \int(\to , \tt)\D_{2}^{*}(\po , \pt)dx.
\]
To verify the Jacobi identity for the bracket induced by $\D_{2}$ we adapt the notation of prolongations (see \cite{Olver} for details). Let $\Theta_{2} $ be the functional bivector associated to $\D_{2}$:
\[ \Theta_{2}= \frac{1}{2}\int (\theta_{1},\theta_{2}) \D_{2}(\theta_{1},\theta_{2})dx
\]
\[ =\frac{1}{2}\int \{ \theta_{1}\wedge\lop\dx\theta_{1}- \theta_{1}\wedge\theta_{2}\dx v \]
\[ +\theta_{2}\wedge\theta_{1}\dx v-\theta_{2}\wedge n\dx \theta_{2}-\theta_{2}\wedge\dx(\theta_{2}n)\} dx
\]
\[ =\frac{1}{2}\int \{ \theta_{1}\wedge\lop\dx\theta_{1}- 2(\dx v)\theta_{1}\wedge\theta_{2}-2\theta_{2}\wedge n\dx \theta_{2}\}dx.
\]
Then $\D_{2}$ is Hamiltonian since
\[ \mbox{pr} \ {\bf v}_{\D_{2}\theta}(\Theta_{2})=-\frac{1}{2}\int \{ 2 (\theta_{1}\dx v-n\dx\theta_{2}-\dx(\theta_{2}n))\wedge\theta_{2}\wedge\dx\theta_{2}
\]
\[ + \dx (\lop\dx\theta_{1}-\theta_{2}\dx v)\wedge\theta_{1}\wedge\theta_{2}\}dx
\]
\[ =-\int \{\lop\dx^{2}\theta_{1}\wedge\theta_{1}\wedge\theta_{2}\} dx
\]
\[ = -\int \{\lop\theta_{1}\wedge\theta_{1}\wedge\theta_{2}\} dx
\]
\[ = -\frac{1}{2}\int \{\dx^{2}(\lop\theta_{1}\wedge\lop\theta_{1})\wedge\theta_{2}\}dx
\]
\[ =0.
\]

These two Hamiltonian structures, (\ref{eq:poi1}) and (\ref{eq:poi2}), are compatible, i.e. the equation (mEP) is bihamiltonian in one space dimension. 
To prove the compatibility it is enough to check that
\begin{equation} 
\mbox{pr} \ {\bf v}_{\D_{1}\theta}(\Theta_{\D_{2}})+ \mbox{pr} \ {\bf v}_{\D_{2}\theta}(\Theta_{\D_{1}})=0
\label{eq:compat}
\end{equation}
holds where $\Theta_{\D_{i}}$ denotes the corresponding bivector for $\D_{i}$. Both summands in (\ref{eq:compat}) vanish:
\[ \mbox{pr} \ {\bf v}_{\D_{1}\theta}(\Theta_{\D_{2}})= \frac{1}{2}\int \{ -2 \dx\theta_{1}\wedge\theta_{2}\wedge\dx\theta_{2}-2 \dx^{2}\theta_{2}\wedge\theta_{1}\wedge\theta_{2}\} dx
\]
\[ = -\int \dx(\theta_{1}\wedge\theta_{2}\wedge\dx\theta_{2}) dx
\]
\[ =0,
\]
and similarly we have
\[ \mbox{pr} \ {\bf v}_{\D_{2}\theta}(\Theta_{\D_{1}})= 0
\]
since
\[ \Theta_{\D_{1}}=\frac{1}{2}\int \{-\theta_{1}\wedge\dx\theta_{2}-\theta_{2}\wedge\dx\theta_{1}\}dx
\]
\[ = -\int \{\theta_{1}\wedge\dx\theta_{2}\}dx.
\]
Therefore the modified Euler-Poisson equation (mEP) is bihamiltonian for $m=1$ and this completes the proof of theorem \ref{th:integrable}.

\hfill $\Box$

\end{document}